\documentclass[final]{siamltex}

\newcommand{\tr}{^{\sf T}}
\newcommand{\m}[1]{{\bf{#1}}}
\newcommand{\g}[1]{\bm #1}
\newcommand{\C}[1]{{\cal {#1}}}

\usepackage{bm}
\usepackage{amsmath}
\usepackage{mathtext}
\usepackage{amsfonts}
\usepackage{amssymb}
\usepackage[mathscr]{eucal}
\usepackage{makeidx}
\usepackage{latexsym}
\usepackage{lscape}
\usepackage{fancyhdr}
\usepackage{graphicx}

\title{An ellipsoidal branch and bound algorithm for global optimization
\thanks{
June 28, 2008.
This material is based upon work supported by the
National Science Foundation under Grants 0619080 and 0620286.} }

\author{
        William W. Hager\thanks{{\tt hager@math.ufl.edu},
        http://www.math.ufl.edu/$\sim$hager,
        PO Box 118105,
        Department of Mathematics,
        University of Florida, Gainesville, FL 32611-8105.
        Phone (352) 392-0281. Fax (352) 392-8357.}
\and
        Dzung T. Phan\thanks{{\tt dphan@math.ufl.edu},
        http://www.math.ufl.edu/$\sim$dphan,
        PO Box 118105,
        Department of Mathematics,
        University of Florida, Gainesville, FL 32611-8105.
        Phone (352) 392-0281. Fax (352) 392-8357.}
}

\begin{document}

\maketitle

\begin{abstract}
A branch and bound algorithm is developed for
global optimization.  Branching in the algorithm
is accomplished by subdividing the feasible set
using ellipses.  Lower bounds are obtained by
replacing the concave part of the objective
function by an affine underestimate.  A ball
approximation algorithm, obtained by generalizing
of a scheme of Lin and Han, is used to solve the
convex relaxation of the original problem.
The ball approximation algorithm is compared
to SEDUMI as well as to gradient projection
algorithms using randomly generated test problems
with a quadratic objective and ellipsoidal constraints.
\end{abstract}



\pagestyle{myheadings} \thispagestyle{plain} \markboth{WILLIAM W.
HAGER AND DZUNG T. PHAN}{ELLIPSOIDAL BRANCH AND BOUND}

\section{Introduction}
\label{introduction}
In this paper we develop a branch and bound algorithm for the
global optimization of the problem
\begin{equation} \tag{P}
\min  \;\; f (\m{x}) \quad
\mbox{subject to } \m{x} \in \Omega ,
\end{equation}
where $\Omega \subset \mathbb{R}^n$ is a compact set and
$f: \mathbb{R}^n \rightarrow \mathbb{R}$ is a weakly convex function
\cite{Vial1983};
that is,
$f(\m{x}) + \sigma \|\m{x}\|^2$ is convex for some $\sigma \ge 0$.
The algorithm starts with a known ellipsoid
$\C{E}$ containing $\Omega$.
The branching process in the branch and
bound algorithm is based on successive ellipsoidal bisections of the
original $\C{E}$. A lower bound for the objective function value
over an ellipse is obtained by writing $f$ as the sum of a convex
and a concave function and replacing the concave part by an affine
underestimate. See \cite{Flu95, HPT95} for discussions concerning
global optimization applications.


As a specific application of our global optimization algorithm,
we consider problems with a quadratic objective function and
with quadratic, ellipsoidal constraints.
Global optimization algorithms for problems with quadratic
objective function and quadratic constraints include those in
\cite{Le00,JLin05,Rab98}.
In \cite{Rab98} Raber considers
problems with nonconvex, quadratic constraints
and with an $n$-simplex enclosing the feasible region.
He develops a branch and bound algorithm based on a
simplicial-subdivision of the feasible set
and a linear programming relaxation over a
simplex to estimate lower bounds.
In a similar setting with box constraints,
Linderoth \cite{JLin05} develops a branch and bound algorithm
in which the the feasible region is subdivided using
the Cartesian product of two-dimensional triangles and rectangles.
Explicit formulae for
the convex and concave envelops of bilinear functions over triangles
and rectangles were derived.
The algorithm of Le \cite{Le00} focuses on problem with convex
quadratic constraints;
Lagrange duality is used to obtain lower bounds for the objective
function, while ellipsoidal bisection is used to subdivide the
feasible region.


The paper is organized as follows. In Section \ref{bisection} we
review the ellipsoidal bisection scheme of \cite{Le00} which is used
to subdivide the feasible region. Section \ref{bounding} develops
the convex underestimator used to obtain a lower bound for the
objective function.
Since $f$ is weakly convex, we can write it as the sum of a convex
and concave functions:
%
\begin{equation}\label{DC}
f(\m{x}) = \left( f (\m{x}) + \sigma \|\m{x}\|^2 \right) +
\left( -\sigma \|\m{x}\|^2 \right) ,
\end{equation}
where $\sigma \ge 0$.
A decomposition of this form is often called a DC (difference
convex) decomposition (see \cite{HPT95}).
For example, if $f$ is a quadratic, then we could take
\[
\sigma = - \min \{ 0, \lambda_1 \} ,
\]
where $\lambda_1$ is the smallest eigenvalue of the Hessian $\nabla^2 f$.
The concave term $-\sigma \|\m{x}\|^2$ in (\ref{DC}) is underestimated
by an affine function $\ell$ which leads to a
convex underestimate $f_L$ of $f$ given by
\begin{equation}\label{fL}
f_L(\m{x}) = \left( f (\m{x}) + \sigma \|\m{x}\|^2 \right)
+ \ell (\m{x}) .
\end{equation}
We minimize $f_L$ over the set $\C{E} \cap \Omega$ to obtain a lower
bound for the objective function on a subset of the feasible set.
An upper bound for the optimal objective function value is obtained
from the best feasible point produced when computing the lower
bound, or from any local algorithm applied to this best feasible
point.
Note that weak convexity for a real-valued function is the analogue of
hypomonotonicity for the derivative operator \cite{cp04, ips03, pen02}.

In Section \ref{PhaseOne} we discuss the phase one problem of
finding a point in $\Omega$ which also lies in the ellipsoid
$\C{E}$. Section \ref{BB} gives the branch and bound algorithm and
proves its convergence. Section \ref{BallReduction} focuses on the
special case where $f$ and $\Omega$ are convex.
The ball approximation algorithm
of Lin and Han \cite{Lin03,LinHan03} for projecting a point onto a
convex set is generalized to replace the norm objective function by
an arbitrary convex function. Numerical experiments, reported in
Section \ref{numerical}, compare the ball approximation algorithm to
SEDUMI 1.1 as well as to gradient projection algorithms. We also
compare the branch and bound algorithm to a scheme of An \cite{Le00}
in which the lower bound is obtained by Lagrange duality.

{\bf Notation.}
Throughout the paper, $\| \cdot \|$ denotes the Euclidian norm.
Given $\m{x}, \m{y} \in \mathbb{R}^n$, $[\m{x}, \m{y}]$ is
the line segment connecting $\m{x}$ and $\m{y}$:
\[
[\m{x}, \m{y}] = \{ (1-t)\m{x} + t\m{y} : 0 \le t \le 1 \} .
\]
The open line segment, which excludes the ends $\m{x}$ and $\m{y}$,
is denoted $(\m{x}, \m{y})$.
The interior of a set $\C{S}$ is denoted $\mbox{int } \C{S}$,
while $\mbox{ri } \C{S}$ is the relative interior.
The gradient $\nabla f (\m{x})$ is a row vector with
\[
(\nabla f (\m{x}))_i = \frac{\partial f (\m{x})}{\partial x_i} .
\]
The diameter of a set $\C{S}$ is denoted $\delta (\C{S})$:
\[
\delta (\C{S}) =
\sup \;\; \{ \|\m{x} - \m{y}\| : \m{x}, \; \; \m{y} \in \C{S} \} .
\]

\section{Ellipsoidal bisection}
\label{bisection}

In this section, we give a brief overview of the ellipsoidal
bisection scheme introduced by An \cite{Le00}.
This idea originates from the ellipsoid method for solving convex optimization
problems by Shor, Nemirovski and Yudin \cite{Shor1977, Yudin1976}.
Consider an ellipsoid $\C{E}$ with center $\m{c}$ in the form
\begin{equation}\label{eq01}
\C{E} = \{\m{x} \in \mathbb{R}^n:
(\m{x}-\m{c})\tr \m{B}^{-1}(\m{x}-\m{c}) \le 1 \} ,
\end{equation}
where $\m{B}$ is a symmetric, positive definite matrix.
Given a nonzero vector $\m{v} \in \mathbb{R}^n$, the sets
\[
H_{-} = \{ \m{x} \in \C{E} : \m{v}\tr\m{x} \le \m{v}\tr\m{c} \}
\quad \mbox{and} \quad
H_{+} = \{ \m{x} \in \C{E} : \m{v}\tr\m{x} \ge \m{v}\tr\m{c} \}
\]
partition $\C{E}$ into two sets of equal volume.
The centers $\m{c}_+$ and $\m{c}_-$ and the matrix
$\m{B}_{\pm}$ of the ellipsoids $\C{E}_{\pm}$
of minimum volume containing $H_{\pm}$ are given as follows:
\[
\m{c}_{\pm} = \m{c} \pm \frac{\m{d}}{n+1}, \quad
\m{B}_{\pm} =
\frac{n^2}{n^2 - 1}\left( \m{B} - \frac{2\m{dd}\tr}{n+1} \right),
\quad \m{d} = \frac{\m{Bv}}{\sqrt{\m{v}\tr\m{Bv}}} .
\]
As mentioned in \cite{Le00}, if the normal $\m{v}$
always points along the major axis of $\C{E}$,
then a nested sequence of bisections shrinks to a point.

\section{Bounding procedure}
\label{bounding}

In this section, we obtain an affine underestimate $\ell$
for the concave function $-\|\m{x}\|^2$ on the
ellipsoid
\begin{equation}\label{ellipsoid}
\C{E} = \{ \m{x} \in \mathbb{R}^n :
\m{x}\tr\m{Ax} - 2\m{b}\tr\m{x} \le \rho \} ,
\end{equation}
where $\m{A}$ is a symmetric, positive definite matrix,
$\m{b} \in \mathbb{R}^n$, and $\rho \in \mathbb{R}$.
The set of affine underestimates for $-\|\m{x}\|^2$
is given by
\begin{equation}\label{under}
\C{U} = \{ \ell: \mathbb{R}^n \rightarrow \mathbb{R}, \;\;
\ell \mbox{ is affine,} \;\; - \|\m{x}\|^2 \ge \ell(\m{x})
\mbox{ for all } \m{x} \in \C{E} \} .
\end{equation}
The best underestimate is a solution of the problem
\begin{equation}\label{linear}
\min_{\ell \in \C{U}} \;\; \max_{\m{x}\in \C{E}} \;\;
- \left( \|\m{x}\|^2 + \ell (\m{x}) \right) .
\end{equation}
\smallskip

\begin{theorem}
\label{UnderTheorem}
A solution of $(\ref{linear})$ is
$\ell^* (\m{x}) = -2\m{c}\tr\m{x} + \gamma$, where
$\m{c} = \m{A}^{-1}\m{b}$ is the center of the ellipsoid,
$\gamma = 2\m{c}\tr\g{\mu} - \|\g{\mu}\|^2$, and
\begin{equation}\label{mu}
\g{\mu} = \mbox{\rm arg } \max_{\m{x} \in \C{E}} \|\m{x} - \m{c}\|^2 .
\end{equation}
If $\delta (\C{E})$ is the diameter of $\C{E}$,
then
\[
\min_{\ell \in \C{U}} \;\; \max_{\m{x}\in \C{E}} \;\;
- \left( \|\m{x}\|^2 + \ell (\m{x}) \right) =
\frac{\delta (\C{E})^2}{4}.
\]
\end{theorem}
\smallskip

\begin{proof}
To begin, we will show that the minimization in (\ref{linear})
can be restricted to a compact set.
Clearly, when carrying out the minimization in (\ref{linear}),
we should restrict our attention to those $\ell$ which touch
the function $h(\m{x}) = -\|\m{x}\|^2$ at some point in $\C{E}$.
Let $\m{y} \in \C{E}$ denote the point of contact.
Since $h(\m{x}) \ge \ell (\m{x})$ and $h (\m{y}) = \ell (\m{y})$,
a lower bound for the error $h (\m{x}) - \ell (\m{x})$ over
$\m{x} \in \C{E}$ is
\[
h (\m{x}) - \ell (\m{x}) \ge |\ell (\m{x}) - \ell (\m{y})|
- |h(\m{x}) - h(\m{y})| .
\]
If $M$ is the
difference between the maximum and minimum value
of $h$ over $\C{E}$, then we have
\begin{equation}\label{change}
h (\m{x}) - \ell (\m{x}) \ge |\ell (\m{x}) - \ell (\m{y})| - M.
\end{equation}

An upper bound for the minimum in (\ref{linear}) is obtained
by the function $\ell_0$ which is constant on $\C{E}$, with value
equal to the minimum of $h(\m{x})$ over $\m{x} \in \C{E}$.
If $\m{w}$ is a point where $h$ attains its minimum over $\C{E}$,
then we have
\[
\max_{\m{x} \in \C{E}} h (\m{x}) - \ell_0 (\m{x}) =
\max_{\m{x} \in \C{E}} h (\m{x}) - h (\m{w}) = M.
\]
For $\m{x} \in \C{E}$, we have
\begin{equation}\label{ell0}
h (\m{x}) - \ell (\m{x}) \le
\max_{\m{x}\in \C{E}} h (\m{x}) - \ell (\m{x}) \le
\max_{\m{x}\in \C{E}} h (\m{x}) - \ell_0 (\m{x}) = M
\end{equation}
when we restrict our attention to affine functions $\ell$ which achieve
an objective function value in (\ref{linear}) which are at least
as good as $\ell_0$.
Combining (\ref{change}) and (\ref{ell0}) gives
\begin{equation}\label{2M}
|\ell (\m{x}) - \ell (\m{y})| \le 2M
\end{equation}
when $\ell$ achieves an objective function value in (\ref{linear})
which is at least as good as $\ell_0$.
Thus, when we carry out the minimization in (\ref{linear}),
we should restrict to affine functions which touch
$h$ at some point $\m{y} \in \C{E}$ and with the change in $\ell$
across $\C{E}$ satisfying the bound (\ref{2M})
for all $\m{x}\in \C{E}$.
This tells us that the minimization in (\ref{linear}) can be
restricted to a compact set, and that a minimizer must exist.

Suppose that $\ell$ attains the minimum in (\ref{linear}).
Let $\m{z}$ be a point in $\C{E}$ where
$h (\m{x}) - \ell (\m{x})$ achieves its maximum.
A Taylor expansion around $\m{x} = \m{z}$ gives
\begin{equation}\label{flz}
h (\m{x}) - \ell (\m{x}) =
h (\m{z}) - \ell (\m{z})  + (\nabla h(\m{z}) - \nabla \ell)(\m{x}-\m{z})
- \|\m{x} - \m{z}\|^2
\end{equation}
since $h(\m{x}) = -\|\m{x}\|^2$.
Since $\ell \in \C{U}$, the set given in (\ref{under}),
we have $h (\m{x}) - \ell (\m{x}) \ge 0$ for all $\m{x} \in \C{E}$,
so (\ref{flz}) yields
\begin{equation}\label{flz2}
0 \le h (\m{z}) - \ell (\m{z})  + (\nabla h(\m{z}) - \nabla \ell)(\m{x}-\m{z})
- \|\m{x} - \m{z}\|^2
\end{equation}
for all $\m{x} \in \C{E}$.
By the first-order optimality conditions for $\m{z}$, we have
\[
(\nabla h(\m{z}) - \nabla \ell)(\m{x}-\m{z}) \le 0
\]
for all $\m{x} \in \C{E}$.
It follows from (\ref{flz2}) that
\[
0 \le h (\m{z}) - \ell (\m{z}) - \|\m{x} - \m{z}\|^2 ,
\]
or
\[
h (\m{z}) - \ell (\m{z})  \ge \|\m{x} - \m{z}\|^2
\]
for all $\m{x} \in \C{E}$.
Since there exists $\m{x} \in \C{E}$ such that
$\|\m{x} - \m{z}\| \ge \delta (\C{E})/2$, we have
\begin{equation}\label{l_lower}
\max_{\m{x} \in \C{E}} h(\m{x}) - \ell (\m{x}) = h (\m{z}) - \ell
(\m{z}) \ge \delta (\C{E})^2/4.
\end{equation}

We now observe that for the specific affine function $\ell^*$
given in the statement of the theorem, (\ref{l_lower}) becomes an
equality, which implies the optimality of $\ell^*$ in (\ref{linear}).
Expand in a Taylor series around $\m{x} = \m{c}$,
where $\m{c} = \m{A}^{-1}\m{b}$ is the center of the ellipsoid $\C{E}$,
to obtain
\[
h(\m{x}) = -\|\m{c}\|^2 - 2\m{c}\tr(\m{x}-\m{c}) -\|\m{x}-\m{c}\|^2
= -2\m{c}\tr\m{x} + \|\m{c}\|^2 - \|\m{x}-\m{c}\|^2.
\]
Hence, for $\ell^*$, we have
\begin{eqnarray*}
h(\m{x}) - \ell^* (\m{x}) &=&
\|\m{c}\|^2 - \gamma - \|\m{x} - \m{c}\|^2 =
\|\g{\mu}-\m{c}\|^2 - \|\m{x}-\m{c}\|^2 \\
&=& \max_{\m{y} \in \C{E}} \|\m{y}-\m{c}\|^2 - \|\m{x}-\m{c}\|^2 .
\end{eqnarray*}
Clearly, $h(\m{x}) - \ell^* (\m{x}) \ge 0$ for all $\m{x} \in \C{E}$,
and the maximum over $\m{x} \in \C{E}$ is attained at $\m{x} = \m{c}$.
Moreover,
\[
h(\m{c}) - \ell^* (\m{c}) =
\max_{\m{y} \in \C{E}} \|\m{y} - \m{c}\|^2 =
\delta (\C{E})^2/4 .
\]
Consequently, (\ref{l_lower}) becomes an equality for $\ell = \ell^*$,
which implies the optimality of $\ell^*$ in (\ref{linear}).
\end{proof}

To evaluate the best affine underestimate given by Theorem
\ref{UnderTheorem}, we need to solve the optimization problem
(\ref{mu}).
This amounts to finding the major axis of the ellipsoid.
The solution is
\[
\g{\mu} = \m{c} + s\m{y},
\]
where $\m{y}$ is a unit eigenvector of $\m{A}$ associated
with the smallest eigenvalue $\epsilon$, and $s$ is chosen so that
$\g{\mu}$ lies on the boundary of the $\C{E}$.
From the definition of $\C{E}$, we obtain
\[
s = \sqrt{(\m{c}\tr\m{Ac} + \rho)/\epsilon} .
\]

We minimize the function $f_L$ in (\ref{fL}) over $\C{E} \cap \Omega$,
with $\ell$ the best affine underestimate of
$-\|\m{x}\|^2$, to obtain a lower bound for the objective
function over $\C{E}\cap \Omega$.
An upper bound for the optimal objective function value
is obtained by starting any local optimization algorithm
from the best iterate generated during the computation of the lower bound.
For the numerical experiments reported later, the
gradient projection algorithm \cite{hz05a} is the local optimization
algorithm.
Of course, by using a faster local algorithm,
the overall speed of the global optimization algorithm will increase.

\section{Phase one}
\label{PhaseOne}
In each step of the branch and bound algorithm for (P),
we need to solve a problem of the form
\begin{equation}\label{LowerBound}
\min \;\; f(\m{x}) \quad
\mbox{subject to } \m{x} \in \C{E} \cap \Omega ,
\end{equation}
in the special case where $f$ is convex (the function
$f_L$ in (\ref{fL})) and
$\C{E}$ is an ellipsoid.
In order to solve this problem,
we often need to find a feasible point.
One approach for finding a feasible point is to
consider the minimization problem
\begin{equation}\label{PhaseOneProb}
\min \;\; \m{x}\tr\m{A}\m{x} - 2 \m{b}\tr\m{x} \quad
\mbox{subject to } \m{x} \in \Omega ,
\end{equation}
where $\m{A}$ and $\m{b}$ are associated with the ellipsoid
$\C{E}$ in (\ref{ellipsoid}).
Assuming we know a feasible point $\m{x}_0 \in \Omega$,
we could apply an optimization algorithm to (\ref{PhaseOneProb}).
If the objective function value can be reduced below $\rho$,
then we obtain a point in $\C{E}$.
If the optimal objective function value is strictly larger
than $\rho$, then the problem (\ref{LowerBound}) is infeasible.

If the set $\Omega$ is itself the intersection of ellipsoids, then
the procedure we have just described could be used in a recursive
fashion to determine a feasible point for either $\Omega$ or $\C{E}
\cap \Omega$, if it exists. In particular, suppose $\Omega =
\cap_{j=1}^m \C{E}_j$ is the intersection of $m$ ellipsoids, where
\[
\C{E}_j = \{ \m{x}\in \mathbb{R}^n : \m{x}\tr\m{A}_j\m{x} - 2
\m{b}_j\tr\m{x} \le \rho_j \} .
\]
A point $\m{x}_1 \in \C{E}_1$ is readily determined.
Proceeding by induction, suppose that we have a point
$\m{x}_{k-1} \in \cap_{j=1}^{k-1} \C{E}_j$.
Any globally convergent iterative method is applied
to the convex optimization problem
\[
\min \;\; \m{x}\tr\m{A}_k\m{x} - 2 \m{b}_k\tr\m{x} \quad
\mbox{subject to } \m{x} \in \cap_{j=1}^{k-1} \C{E}_j .
\]
If the objective function value is reduced below $\rho_k$,
then a feasible point in $\cap_{j=1}^k \C{E}_j$ has been determined.
Conversely, if the optimal objective function value is above
$\rho_k$, then $\cap_{j=1}^k \C{E}_j$ is empty.

\section{Branch and bound algorithm}
\label{BB}
Our branch and bound algorithm is patterned after a
general branch and bound algorithm, as appears in \cite{HPT95}
for example.
For any ellipse $\C{E}$, define
\begin{equation}\label{l}
M_L (\C{E}) =
\min \;\; \{ f_L (\m{x}) : \m{x} \in \C{E} \cap \Omega \} ,
\end{equation}
where $f_L$ is the lower bound (\ref{fL}) corresponding to the best
affine underestimate of $-\|\m{x}\|^2$ on $\C{E}$.
We assume that an algorithm is available to solve the
optimization problem (\ref{l}).
\smallskip

\begin{itemize}
\item[] \hspace{-.2in}\textbf{Ellipsoidal branch and bound with linear
underestimate (EBL)}

\item [1.]
Let $\C{E}_0$ be an ellipsoid which contains $\Omega$
and set $\C{S}_0 = \{\C{E}_0\}$.

\item [2.]
Evaluate $M_L (\C{E}_0)$ and let $\m{x}_0 \in \Omega$ denote the feasible
point generated during the evaluation of $M_L(\C{E}_0)$ with
the smallest function value.

\item [3.]
For $k = 0, 1, 2, \ldots $
\begin{itemize}
\item [(a)]
Choose $\C{E}_k \in \C{S}_k$ such that
$M_L (\C{E}_k) =  \min \{ M_L (\C{E}) : \C{E} \in \C{S} \}$.
Bisect $\C{E}_k \in \C{S}_k$ with
two ellipsoids denoted $\C{E}_{k1}$ and $\C{E}_{k2}$
(see Section \ref{bisection}).
Evaluate $M_L (\C{E}_{k1})$ and $M_L (\C{E}_{k2})$.
\item [(b)]
Let $\m{x}_{k+1}$ denote a feasible point associated with the
smallest function value that has
been generated up to this iteration and up to this step.
Hence, if $\m{y}_{k1}$ and $\m{y}_{k2}$ are solutions to
(\ref{l}) associated with $\C{E} = \C{E}_{k1}$ and $\C{E} = \C{E}_{k2}$
respectively, then we have $f(\m{x}_{k+1}) \le f(\m{y}_{ki})$,
$i = 1, 2$.
\item [(c)]
Set $\C{S}_{k+1} = \{ \C{E} \in \C{S}_k \cup
\{ \C{E}_{k1} \} \cup \{ \C{E}_{k2}\}  : M_L (\C{E}) \le f(\m{x}_{k+1}),
\C{E} \ne \C{E}_k \}$
\end{itemize}
\end{itemize}
\bigskip

\begin{theorem} \label{theo3}
Suppose that the following conditions hold:
\begin{itemize}
\item [A1.]
The feasible set $\Omega$ is contained in some given ellipsoid
$\C{E}$, $\Omega$ is compact, and $f$ is weakly convex
over $\C{E}$.
\item [A2.]
A nested sequence of ellipsoidal bisections shrinks to a point
(see Section \ref{bisection}).
\end{itemize}
Then
every accumulation point
of the sequence $\m{x}_k$ is a solution of {\rm (P)}.
\end{theorem}
\smallskip

\begin{proof}
Let $\m{y}$ denote any global minimizer for (P).
We now show that for each $k$,
there exists $\C{E} \in \C{S}_k$ with $\m{y} \in \C{E}$.
Since $\Omega \subset  \C{E}_0$, $\m{y} \in \C{E}_0$.
Proceeding by induction, suppose that for each $j$,
$0 \le j \le k$, there exists an ellipsoid
$\C{F}_j \in \C{S}_j$ with $\m{y} \in \C{F}_j$.
We now wish to show that there exist $\C{F}_{k+1} \in \C{S}_{k+1}$
with $\m{y} \in \C{F}_{k+1}$.
In Step 3c, $\C{F}_k \in \C{S}_k$
can only be deleted from $\C{S}_{k+1}$ if
$M_L(\C{F}_k) > f(\m{x}_{k+1})$ or $\C{F}_k = \C{E}_k$.
The former case cannot occur since
\[
M_L(\C{F}_k) \le f(\m{y}) \le f(\m{x}_{k+1}),
\]
due to the global optimality of $\m{y}$.
If $\C{F}_k = \C{E}_k$, then
$\m{y}$ lies in either $\C{E}_{k1}$ or $\C{E}_{k2}$.
If $\m{y} \in \C{E}_{ki}$, then $\C{E}_{ki} \in \C{S}_{k+1}$ since
\[
M_L(\C{E}_{ki}) \le f(\m{y}) \le f(\m{x}_{k+1}).
\]

Let $\m{x}^*$ denote
an accumulation point of the sequence $\m{x}_k$.
Since $\Omega$ is closed and $\m{x}_k \in \Omega$ for each $k$,
$\m{x}^* \in \Omega$.
By \cite[Prop. 4.4]{Vial1983}, a weakly convex function is locally
Lipschitz continuous.
Hence, $f$ is continuous on $\Omega$ and
$f (\m{x}_k)$ approaches $f (\m{x}^*)$.
If $\m{x}^*$ is a solution of (P), then the proof is complete.
Otherwise, $f (\m{y}) < f (\m{x}^*)$.

For each $k$, we have
\begin{equation}\label{EL}
\min \;\; \{ M_L (\C{E}) : \C{E} \in \C{S}_k \} \le
M_L(\C{F}_k) \le f (\m{y}) < f(\m{x}^*) .
\end{equation}
Let $\C{G}_k$ denote an ellipsoid which achieves the minimum on the
left side of (\ref{EL}) and let $\m{y}_k$ denote a minimizer in
(\ref{l}) corresponding to $\C{E} = \C{G}_k$.
The inequality (\ref{EL}) reduces to
\begin{equation}\label{L1}
f_L(\m{y}_k) \le f (\m{y}) < f(\m{x}^*) .
\end{equation}
Since $\m{y}$ minimizes $f$ over $\Omega$, (\ref{L1}) implies that
\begin{equation}\label{L2}
f_L (\m{y}_k) \le f (\m{y}) \le f (\m{y}_k) .
\end{equation}
By Theorem \ref{UnderTheorem},
\begin{equation}\label{L3}
f(\m{y}_{k}) - f_L(\m{y}_{k}) =
-\sigma (\|\m{y}_k\|^2 + \ell (\m{y}_k))
\le \sigma \delta (\C{G}_{k})^2/4 ,
\end{equation}
where $\ell$ is the best linear lower bound for the function
$-\|\m{x}\|^2$, and $\sigma \ge 0$ is the parameter associated
with the convex/concave decomposition (\ref{DC}).

Each ellipsoid $\C{E}_k$ corresponds to a vertex on the
branch and bound tree associated with EBL.
Choose the iteration numbers $k_1 < k_2 < \ldots$
so that they correspond to vertices along an infinite path
on the branch and bound tree, starting from the root of the tree.
By (A2), $\delta (\C{G}_{{k}_i})$ tends to 0 as $i$ tends to infinity.
Hence, (\ref{L3}) implies that
$|f(\m{y}_{{k}_i}) - f_L(\m{y}_{{k}_i})|$ tends to zero.
Combining this with (\ref{L1}) and (\ref{L2}) shows that
$f(\m{y}_{{k}_i}) < f (\m{x}^*)$ for $i$ sufficiently large,
which violates Step 3b and the fact that $f(\m{x}_{k+1})$ is the
smallest function value at step $k$ and the smallest values
monotonically approach $f(\m{x}^*)$.
\end{proof}

Note that if for any $k$,
$f(\m{x}_k) = \min \{ M_L (\C{E}) : \C{E} \in \C{S}_k \}$,
then $\m{x}_k$ is a global minimizer.
\section{Ball approximation algorithm for convex optimization}
\label{BallReduction}

In this section we give an algorithm to solve (P) in the special case that
$f$ and $\Omega$ are convex.
This algorithm, which is based on the successive approximation of the
feasible set by balls, ties in nicely with the ellipsoidal-based
branch and bound algorithm.
The algorithm is a generalization of the ball approximation
algorithm \cite{LinHan03} of Lin and Han.
The algorithm of Lin and Han deals with
the special case where the objective function has the form
$\|\m{x} - \m{a}\|^2$ and $\Omega$ is an intersection of ellipsoids.
Lin generalizes this algorithm in \cite{Lin03} to treat
convex constraints.
The analysis in \cite{Lin03,LinHan03} is tightly
coupled to the norm objective function.
In our further generalization of the Lin/Han algorithm,
the norm objective function is replaced by an arbitrary convex
functional $f$ and an additional constraint set $\chi \subset \mathbb{R}^n$
is included, which might represent bound constraints for example.
More precisely, we consider the problem
\begin{equation} \tag{C}
\min f (\m{x}) \quad
\mbox{subject to }
\m{x} \in \C{F} := \{ \m{x}\in \chi : \m{g}(\m{x}) \le \m{0} \},
\end{equation}
where $f: \mathbb{R}^n \rightarrow \mathbb{R}$,
$\m{g} : \mathbb{R}^n \rightarrow \mathbb{R}^m$, and the
following conditions hold:
\begin{itemize}
\item[{\rm C1.}]
$f$ and $\m{g}$ are convex and differentiable, $\chi$ is closed
and convex, and $\C{F}$ is compact.
\item[{\rm C2.}]
There exists $\bar{\m{x}}$ in the relative
interior of $\chi$ with $\m{g}(\bar{\m{x}}) < \m{0}$.
\item[{\rm C3.}]
There exists $\gamma > 0$ such that
$\|\nabla g_i (\m{x})\| \ge \gamma$ when $g_i (\m{x}) = 0$ for some
$i \in [1, m]$ and $\m{x} \in \chi$.
\end{itemize}
The condition C2 is referred to as the Slater condition.

We will give a new analysis which handles this more general convex problem (C).
In each iteration of Lin's algorithm in \cite{Lin03},
the convex constraints are approximated by ball constraints.
Let $h: \mathbb{R}^n \rightarrow \mathbb{R}$ be a convex,
differentiable function which defines a convex,
nonempty set
\[
\C{H} = \{ \m{x} \in \mathbb{R}^n : h (\m{x}) \le 0 \} .
\]
The ball approximation $\C{B}_h (\m{x})$
at $\m{x} \in \C{H}$ is expressed in
terms of a center map $\m{c}: \C{H} \rightarrow \mathbb{R}^n$
and a radius map $r: \C{H} \rightarrow \mathbb{R}$:
\[
\C{B}_h (\m{x}) = \{ \m{y} \in \mathbb{R}^n :
\| \m{y} - \m{c} (\m{x}) \| \le r (\m{x}) \} .
\]
These two maps must satisfy the following conditions:
\smallskip

\begin{itemize}
\item[\rm B1.] Both $\m{c}$ and $r$ are continuous on $\C{H}$.
\item[\rm B2.] If $h (\m{x}) < 0$,
then $\m{x} \in \mbox{int } \C{B}_h (\m{x})$, the interior of
$\C{B}_h (\m{x})$.
\item[\rm B3.] If $h (\m{x}) = 0$,
then $\m{x} \in \partial \C{B}_h (\m{x})$,
and $c (\m{x}) = \m{x} - \alpha \nabla h (\m{x})\tr$
for some fixed $\alpha > 0.$
\end{itemize}
\smallskip
Maps which satisfy
B1, B2, and B3 are the following, assuming $h$ is continuously
differentiable:
\[
\m{c}(\m{x}) = \m{x} - \alpha \nabla h (\m{x})\tr, \quad
r (\m{x}) = \alpha \|\nabla h (\m{x})\| - \beta h (\m{x}),
\]
where $\alpha$ and $\beta$ are fixed positive scalars.

Let $\m{c}_i$ and $r_i$ denote center and radius maps associated
with $g_i$, let $\C{B}_i$ be the associated ball given by
\[
\C{B}_i (\m{x}) = \{ \m{y} \in \mathbb{R}^n :
\| \m{y} - \m{c}_i (\m{x}) \| \le r_i (\m{x}) \} ,
\]
and define
$\C{B}(\m{x}) = \cap_{i=1}^m \C{B}_i (\m{x})$.
Our generalization of the algorithm of Lin and Han is the following:
\smallskip

\begin{enumerate}
\item[] \hspace{-.2in}\textbf{Ball approximation algorithm (BAA)}
\item
Let $\m{x}_0$ be a feasible point for (C).

\item For $k = 0, 1, \ldots $
\begin{itemize}
\item[{\rm (a)}]
Let $\m{y}_k$ be a solution of the problem
\begin{equation} \label{eqa2}
\min  \;\; f (\m{x})
\quad \mbox{subject to } \m{x} \in \chi \cap \C{B} (\m{x}_k) .
\end{equation}
\item[{\rm (b)}]
Set $\m{x}_{k + 1} = \m{x}(\tau_k)$ where
$\m{x}(\tau) = (1-\tau)\m{x}_k + \tau \m{y}_k$
and $\tau_k$ is the largest $\tau \in [0, 1]$ such that
$\m{x}(\sigma) \in \C{F}$
for all $\sigma \in [0, \tau]$.
\end{itemize}
\end{enumerate}
\bigskip

In \cite[Lem. 3.1]{Lin03} it is shown that
$\mbox{int } \C{B}(\m{x}) \ne \emptyset$ for each $\m{x} \in \C{F}$
when the center and radius
maps $\m{c}_i$ and $r_i$ satisfy B2 and B3 and there exists
$\bar{\m{x}}$ such that $\m{g}(\bar{\m{x}}) < \m{0}$.
Lin's proof is based on the following observation: For $\epsilon > 0$
sufficiently small, $\m{x} + \epsilon (\bar{\m{x}} - \m{x})$ lies
in the interior of $\C{B}_i (\m{x})$ for each $i$.
In C2 we also assume that $\bar{\m{x}} \in \mbox{ri } \chi$,
where ``ri'' denotes relative interior.
Hence, for $\epsilon > 0$
sufficiently small, $\m{x} + \epsilon (\bar{\m{x}} - \m{x})$ lies
in both $\mbox{ri } \chi$ and in the
interior of $\C{B}_i (\m{x})$ for each $i$.
Consequently, we have
\begin{equation}\label{Slater}
\mbox{ri } \chi \cap \mbox{int } \C{B}(\m{x}) \ne \emptyset \mbox{ for every }
\m{x} \in \C{F} .
\end{equation}
This implies that the subproblems (\ref{eqa2}) of BAA are always feasible.
An optimal solution $\m{y}_k$ exists due to the compactness of the feasible
set and the continuity of the objective function.
\smallskip

\begin{theorem}
If {\rm C1}, {\rm C2}, and {\rm C3} hold and
the center map $\m{c}_i$ and the radius map $r_i$ satisfy
{\rm B1}, {\rm B2}, and {\rm B3}, $i = 1, 2, \ldots , m$,
then the limit $\m{x}^*$ of any
convergent subsequence of iterates
$\m{x}_k$ of Algorithm $2$ is a solution of {\rm (C)}.
\end{theorem}
\smallskip

\begin{proof}
Initially, $\m{x}_0 \in \C{F}$.
Proceeding by induction, it follows from the line search in
Step 2a of BAA that $\m{x}_k \in \C{F}$ for each $k$.
By B2 and B3,
$\m{x} \in \C{B}_{h} (\m{x})$ if $h (\m{x}) \le 0$.
Consequently, $\m{x}_k \in \chi \cap \C{B}(\m{x}_k)$ for each $k$.
This shows that $\m{x}_k$ is feasible in (\ref{eqa2}) for each $k$,
and the minimizer $\m{y}_k$ in (\ref{eqa2}) satisfies
\begin{equation}\label{descent}
f (\m{y}_k) \le f (\m{x}_k) \quad \mbox{for each } k.
\end{equation}
By the convexity of $f$ and by (\ref{descent}), we have
\begin{equation}\label{monotone}
f (\m{x}_{k+1})
\le \tau_k f (\m{y}_k) + (1-\tau_k) f (\m{x}_k)
\le f(\m{x}_k) ,
\end{equation}
where $\tau_k \in [0, 1]$ is defined in Step 2b of BAA.
Hence, $f (\m{x}_k)$ approaches a limit monotonically.
Since $\C{F}$ is compact and $\m{x}_k \in \C{F}$ for each $k$,
an accumulation point $\m{x}^* \in \C{F}$ exists.
Since the center maps $\m{c}_i$ and the radius maps $r_i$ are
continuous, the balls $\C{B}_i (\m{x}_k)$ are uniformly bounded,
and hence, the $\m{y}_k$ are contained in bounded set.
Let $\m{y}^*$ denote an accumulation point of the $\m{y}_k$.
To simplify the exposition, let $(\m{x}_k, \m{y}_k)$ denote a pruned
version of the original
sequence which approaches the limit $(\m{x}^*, \m{y}^*)$.

We now show that
\begin{equation}\label{y*}
\m{y}^* = \mbox{arg } \min \;\; f (\m{x})
\quad \mbox{subject to } \m{x} \in \chi \cap \C{B} (\m{x}^*) .
\end{equation}
Suppose, to the contrary, that there exists
$\tilde{\m{y}} \in \chi \cap \C{B}(\m{x}^*)$ such that
$f(\tilde{\m{y}}) < f (\m{y}^*)$.
Referring to the discussion before (\ref{Slater}), choose
$\tilde{\m{x}} \in \mbox{ri } \chi$ with
$\tilde{\m{x}} \in \mbox{int } \C{B}_i (\m{x}^*)$ for each $i$.
Define $\hat{\m{y}} = \tilde{\m{y}} + \epsilon (\tilde{\m{x}} - \tilde {\m{y}})$
where $\epsilon > 0$ is small enough that
$\hat{\m{y}} \in \mbox{ri } \chi$,
$\hat{\m{y}} \in \mbox{int } \C{B}_i (\m{x}^*)$ for each $i$, and
$f(\hat{\m{y}}) < f (\m{y}^*)$.
For $k$ sufficiently large,
$\hat{\m{y}} \in \chi \cap \C{B}(\m{x}_k)$ due to the
continuity of the center and radius maps.
Since $f (\m{y}_k)$ approaches $f (\m{y}^*) > f (\tilde{\m{y}})$,
we contradict the optimality of $\m{y}_k$ in (\ref{eqa2}).
This establishes (\ref{y*}).

Again, by B2 and B3, $\m{x}^*$ is feasible in (\ref{y*}).
Since $\m{y}^*$ is optimal in
(\ref{y*}), we have $f (\m{y}^*) \le f (\m{x}^*)$.
We will show that
\begin{equation}\label{equality}
f (\m{y}^*) = f (\m{x}^*).
\end{equation}
Suppose, to the contrary, that $f (\m{y}^*) < f (\m{x}^*)$.
Since $\m{x}^* \in \C{F}$, we conclude that for each $i$, one
of the following two cases can occur:
\begin{itemize}
\item [(i)]
$g_i (\m{x}^*) = 0$: In this case, it follows from B3 that
$\m{x}^* \in \partial \C{B}_i (\m{x}^*)$.
Since both $\m{x}^*$ and $\m{y}^* \in \chi \cap \C{B}_i (\m{x}^*)$,
we have $[\m{x}^*, \m{y}^*] \in \chi \cap \C{B}_i (\m{x}^*)$.
Hence, the vector $\m{y}^* - \m{x}^*$ makes an acute angle with
the inward pointing normal at $\m{x}^*$.
By B3 the inward pointing normal is a positive multiple of
$-\nabla g (\m{x}^*)$; it follows that
\[
-\nabla g (\m{x}^*)(\m{y}^* - \m{x}^*) > 0 .
\]
By a Taylor expansion around $\m{x}^*$, we see that
there exist $\sigma_i \in (0, 1)$ such that
\begin{equation}\label{taui}
g_i (\m{x}^* + \sigma(\m{y}^*-\m{x}^*)) < 0 \quad
\mbox{for all } \sigma \in (0, \sigma_i] .
\end{equation}
\item[(ii)]
$g_i (\m{x}^*) < 0$: In this case, there trivially exists
$\sigma_i \in (0, 1)$ such that (\ref{taui}) holds.
\end{itemize}

Let $\sigma^*$ be the minimum of $\sigma_i$, $1 \le i \le m$.
By the convexity of $f$, we have
\begin{equation}\label{strict}
f (\m{x}^* + \sigma^* (\m{y}^* - \m{x}^*)) \le
f (\m{x}^*) + \sigma^* (f(\m{y}^*) - f(\m{x}^*)) < f (\m{x}^*)
\end{equation}
since $f(\m{y}^*) < f(\m{x}^*))$.
Since both $\m{x}_k$ and $\m{y}_k \in \chi \cap \C{B}(\m{x}_k)$,
the line segment $[\m{x}_k, \m{y}_k]$ is contained in
$\chi \cap \C{B}(\m{x}_k)$.
Since $\m{x}_k$ approaches $\m{x}^*$ and $\m{y}_k$ approaches $\m{y}^*$,
it follows from (\ref{taui}) that
\[
\m{x}_k + \sigma^* (\m{y}_k - \m{x}_k) \in \C{F}
\]
for $k$ sufficiently large.
Again, by the convexity of $f$, (\ref{descent}),
and the fact that $\tau_k$ is taken
as large as possible so that
\[
\m{x}_{k+1} = (1-\tau_k)\m{x}_k + \tau_k \m{y}_k \in \C{F} ,
\]
we have
\begin{equation}\label{h1}
f (\m{x}_{k+1}) \le f(\m{x}_k) + \tau_k (f(\m{y}_k) - f (\m{x}_k))
\le f(\m{x}_k) + \sigma^* (f(\m{y}_k) - f (\m{x}_k)) .
\end{equation}
Since $(\m{x}_k, \m{y}_k)$ converges to $(\m{x}^*, \m{y}^*)$,
it follows from (\ref{strict}) that
\begin{equation}\label{h2}
\lim_{k \rightarrow \infty}
f(\m{x}_k) + \sigma^* (f(\m{y}_k) - f (\m{x}_k)) =
f(\m{x}^*) + \sigma^* (f(\m{y}^*) - f (\m{x}^*)) < f (\m{x}^*) .
\end{equation}
Hence, for $k$ sufficiently large, (\ref{h1}) and (\ref{h2}) imply that
$f (\m{x}_{k+1}) < f (\m{x}^*)$,
which contradicts the monotone decreasing convergence
(\ref{monotone}) of $f (\m{x}_k)$ to $f(\m{x}^*)$.
This completes the proof of (\ref{equality}).

Let $L : \mathbb{R}^{m+n} \rightarrow \mathbb{R}$ be
the Lagrangian defined by
\[
L (\g{\lambda}, \m{x}) =
f (\m{x}) +
\frac{1}{2}\sum_{i=1}^m
\lambda_i \left( \|\m{x} - \m{c}_i (\m{x}^*)\|^2 - r_i (\m{x}^*)
\right) .
\]
Since $\m{y}^*$ is a solution of (\ref{y*}) and the Slater condition
(\ref{Slater}) holds,
the first-order optimality condition holds at $\m{y}^*$.
That is, there exist $\g{\lambda}^* \in \mathbb{R}^m$ such that
\begin{equation}\label{1st-order}
\left.
\begin{array}{l}
\g{\lambda}^* \ge \m{0},
\quad \lambda_i^* (\|\m{y}^* - \m{c}_i (\m{x}^*)\|^2 - r_i^2 (\m{x}^*)) = 0,
\;\; i = 1, 2, \ldots , m, \\[.1in]
\nabla_x L (\g{\lambda}^*, \m{y}^*)(\m{x} - \m{y}^*) \ge \m{0}
\mbox{ for all } \m{x} \in \chi .
\end{array}
\right\}
\end{equation}
If $\nabla f (\m{y}^*)(\m{x}-\m{y}^*) \ge \m{0}$ for all $\m{x} \in \chi$,
then $\m{y}^*$ is the global minimizer of the convex function $f$ over
$\chi$.
Since $f (\m{y}^*) = f (\m{x}^*)$ by (\ref{equality}),
it follows that $\m{x}^*$
is a solution of (C), and the proof would be complete.
Hence, we suppose that $\nabla f (\m{y}^*)(\m{x}-\m{y}^*) < \m{0}$
for some $\m{x} \in \chi$,
which implies that $\g{\lambda}^* \ne \m{0}$ by (\ref{1st-order}).

Since $f$ is convex, we have
\begin{equation}\label{c1}
f (\m{x}^*) \ge f (\m{y}^*) + \nabla f (\m{y}^*)(\m{x}^* - \m{y}^*) .
\end{equation}
We expand the expression
\[
\frac{1}{2}\sum_{i=1}^m
\lambda_i \left( \|\m{x} - \m{c}_i (\m{x}^*)\|^2 - r_i^2 (\m{x}^*) \right)
\]
in a Taylor series around $\m{x} = \m{y}^*$ and evaluate at
$\m{x} = \m{x}^*$ to obtain
\begin{eqnarray*}
\frac{1}{2}\sum_{i=1}^m
\lambda_i \left( \|\m{x}^* - \m{c}_i (\m{x}^*)\|^2 - r_i^2 (\m{x}^*) \right)
&=&
\frac{1}{2}\sum_{i=1}^m
\lambda_i \left( \|\m{y}^* - \m{c}_i (\m{x}^*)\|^2 - r_i^2 (\m{x}^*) \right) \\
&& \hspace {-1in} + \sum_{i=1}^m \lambda_i (\m{y}^* - \m{c}_i (\m{x}^*))\tr
(\m{x}^* - \m{y}^*) +
\frac{1}{2} \|\m{x}^* - \m{y}^*\|^2 \sum_{i=1}^m \lambda_i^* .
\end{eqnarray*}
We add this equation to (\ref{c1}) to obtain
\begin{equation}\label{L_lower}
L(\m{x}^*, \g{\lambda}^*) \ge
L(\m{y}^*, \g{\lambda}^*) +
\nabla_x L(\m{y}^*, \g{\lambda}^*)(\m{x}^* - \m{y}^*) +
\frac{1}{2} \|\m{x}^* - \m{y}^*\|^2 \sum_{i=1}^m \lambda_i^* .
\end{equation}
By complementary slackness and by (\ref{equality}),
we have
$L(\m{y}^*, \g{\lambda}^*) = f (\m{y}^*) = f (\m{x}^*)$.
Hence, (\ref{L_lower}) yields
\begin{eqnarray}
\frac{1}{2}
\|\m{x}^* - \m{y}^* \|^2
\sum_{i=1}^m \lambda_i^* &\le&
-\nabla_x L(\m{y}^*, \g{\lambda}^*)(\m{x}^* - \m{y}^*) \nonumber \\
&& \quad +
\frac{1}{2} \sum_{i=1}^m \lambda_i^*
\left( \|\m{x}^* - \m{c}_i (\m{x}^*)\|^2 - r_i^2 (\m{x}^*) \right) .
\label{xy}
\end{eqnarray}
By (\ref{1st-order}) and the fact that $\m{x}^* \in \chi$,
we have $\nabla_x L(\m{y}^*, \g{\lambda}^*)(\m{x}^* - \m{y}^*) \ge \m{0}$.
Since $\m{x}^* \in \C{B}(\m{x}^*)$, the last term in (\ref{xy}) is
nonpositive.
Hence, the entire right side of (\ref{xy}) is nonpositive.
Since $\g{\lambda}^* \ge \m{0}$
and $\g{\lambda^*} \ne \m{0}$,
(\ref{xy}) implies that $\m{y}^* = \m{x}^*$.

Replacing $\m{y}^*$ by $\m{x}^*$ in the first-order
conditions (\ref{1st-order}) gives
\begin{equation}\label{B}
\left( \nabla f (\m{x}^*)
+ \sum_{i=1}^m \lambda_i^* (\m{x}^* - \m{c}_i (\m{x}^*))\tr \right)
(\m{x} - \m{x}^*) \ge \m{0} \mbox{ for all } \m{x} \in \chi .
\end{equation}
If $g_i (\m{x}^*) < 0$,
then by B2, $\m{x}^* \in \mbox{ int } \C{B}_i (\m{x}^*)$ and
$\lambda_i^* = 0$ by complementary slackness.
If $g_i(\m{x}^*) = 0$, then by B3,
$\m{c}_i (\m{x}^*) = \m{x}^* - \alpha \nabla g_i (\m{x}^*)\tr$.
With these substitutions, (\ref{B}) yields
\[
\left( \nabla f (\m{x}^*)
+ \alpha \sum_{i=1}^m \lambda_i^* \nabla g_i (\m{x}^*) \right)
(\m{x} - \m{x}^*) \ge \m{0} \mbox{ for all } \m{x} \in \chi .
\]
Hence, the first-order optimality conditions
for (C) are satisfied at $\m{x}^*$.
Since the objective function and the constraints of (C) are
convex, $\m{x}^*$ is a solution of (C).
This completes the proof.
\end{proof}


\section{Numerical experiments}
\label{numerical}
We investigate the performance of the algorithms of the previous
sections using randomly generated quadratically constrained
quadratic programming problems of the form
\begin{equation} \tag{QP}
 \min  \m{x}\tr\m{A}_0 \m{x} + \m{b}_0\tr \m{x}
\mbox { subject to } \m{g}(\m{x}) \le \m{0},
 \end{equation}
where $\m{x}\in \mathbb{R}^n$ and
$g_i (\m{x}) = \m{x}\tr\m{A}_i \m{x} + \m{b}_i\tr \m{x} + c_i$,
$i = 1, 2, \ldots , m$.
Here $\m{b}_i \in \mathbb{R}^n$ and $c_i$ is a scalar for each $i$.
The matrices $\m{A}_i$ are symmetric, positive definite for $i \ge 1$.
In our experiments with the ball approximation algorithm,
we take $\m{A}_0$ symmetric, positive semidefinite.
In our experiments with the branch and bound algorithm, we consider
more general indefinite $\m{A}_0$.
The codes are written in either C or Fortran.
The experiments were implemented using a Matlab 7.0.1
interface on a PC with 2GB memory and Intel Core 2 Duo 2Ghz
processors running the Windows Vista operating system.

\subsection{Rate of convergence for BAA}
\begin{figure}[htp]
\centering
\includegraphics[scale=0.75]{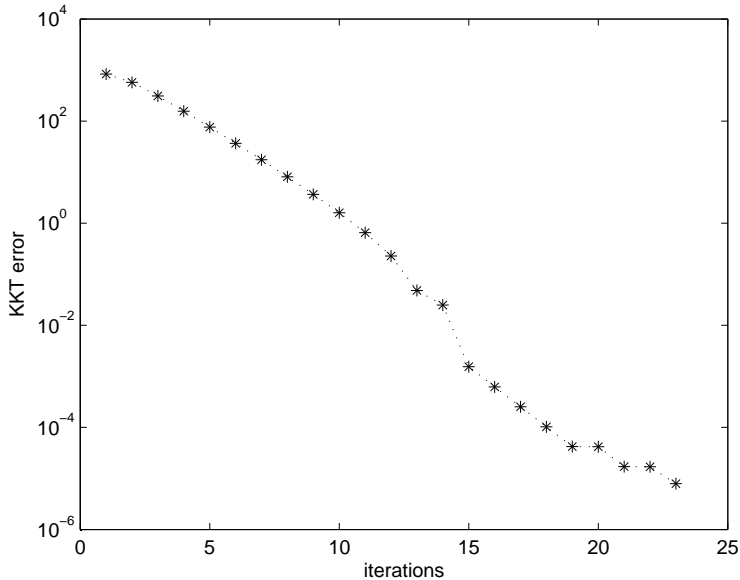}
\caption{KKT error versus iteration number for n = 200, m = 4, and
$\m{A}_0$ positive definite}\label{figpd}
\end{figure}

\begin{figure}[htp]
\centering
\includegraphics[scale=0.75]{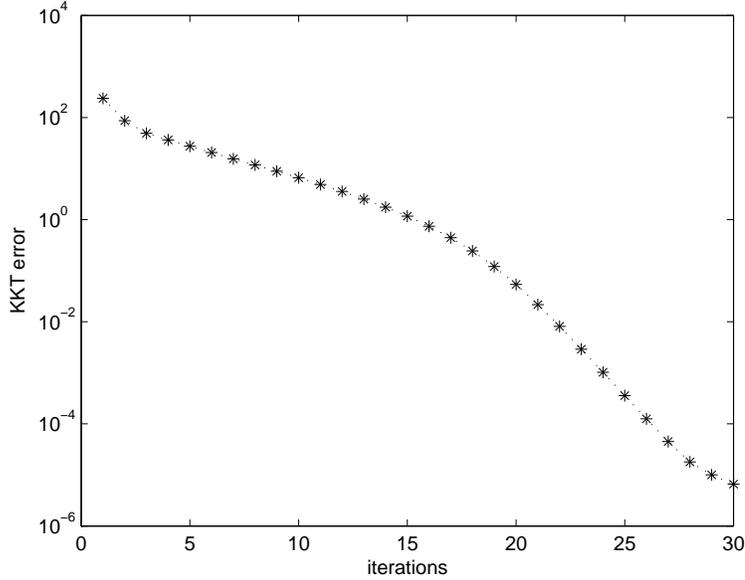}
\caption{KKT error versus iteration number for n = 200, m = 4, and
a positive semidefinite $\m{A}_0$.}\label{figpsd}
\end{figure}

The theory of Section \ref{BallReduction} establishes the
convergence of BAA.
Experimentally, we observe that the convergence rate is linear.
Figure \ref{figpd} shows that the behavior of the KKT error
as a function of the iteration number for a randomly generated
positive definite matrix $\m{A}_0$ of dimension 200 and
for 4 ellipsoidal constraints ($m = 4$).
The KKT error is computed using the formula given in
Section 4 of \cite{HagerGowda99}.
Roughly, this formula amounts to the
infinity norm of the gradient of the Lagrangian
plus the infinity norm of the violation in complementary slackness.
If $\m{A}_0$ is constructed to have
precisely one zero eigenvalue, then the convergence rate again
appears to be linear, as seen in Figure \ref{figpsd}.

\subsection{Comparison with other algorithms for programs with convex cost}
To gain some insight into the relative performance of the ball approximation
algorithm (BAA), we solved randomly generated problems with convex
cost using three other algorithms:
\begin{itemize}
\item
SEDUMI, for optimization over symmetric cones.
\item
The gradient projection algorithm.
We tried both the nonmonotone gradient project algorithm (NGPA)
given in \cite{hz05a} and the nonmonotone spectral projected gradient
method (SPG) of Birgin, Mart\'{\i}nez, and Raydan
\cite{bmr00,bmr01} (ACM Algorithm 813).
\end{itemize}

We now discuss in detail how each of these algorithms was implemented.
The BAA subproblems (\ref{eqa2}) have the form
\begin{equation}\label{BAA}
\min \m{x}\tr \m{A}_0\m{x} + \m{b}_0\tr\m{x} \quad
\mbox{ subject to }
\|\m{x}-\m{c}_i\|^2 \le r_i^2, \quad 1 \le i \le m .
\end{equation}
We solve these subproblems by applying the active set algorithm
(ASA) developed in \cite{hz05a} to the dual problem.
To facilitate the evaluation of the dual function,
we compute the diagonalization $\m{A}_0 = \m{QDQ}\tr$ where $\m{D}$
is diagonal and $\m{Q}$ is orthogonal.
Substituting $\m{x} = \m{Qy}$ in (\ref{BAA}) yields the equivalent
problem
\[
\min \m{y}\tr \m{D}\m{y} + \m{b}_0\tr\m{Qy} \quad
\mbox{ subject to }
\|\m{y}-\m{Q}\tr\m{c}_i\|^2 \le r_i^2, \quad 1 \le i \le m .
\]
The dual problem is
\begin{equation}\label{dual}
\max_{\g{\lambda} \ge \m{0}} \;\; \min_{\m{y}\in \mathbb{R}^n} \;\;
\m{y}\tr\m{Dy} + \m{b}_0\tr\m{Qy} + \sum_{i=1}^m \lambda_i
\left( \|\m{y} - \m{Q}\tr\m{c}_i\|^2 - r_i^2 \right) .
\end{equation}
The $i$-th component of the gradient of the dual function
with respect to $\g{\lambda}$ is simply
$\|\m{y}(\g{\lambda}) - \m{Q}\tr \m{c}_i\|^2 - r_i^2$ where
$\m{y}(\g{\lambda})$ achieves the minimum in (\ref{dual}).
This minimum is easily evaluated since the quadratic term in the
objective function is diagonal.

SEDUMI could be applied directly to (QP) when the cost function
is strongly convex.
We used Version 1.1 of the code obtained from
\smallskip
\begin{center} http://sedumi.mcmaster.ca/ \end{center}
\smallskip
In implementing the gradient projection algorithm for (QP),
we need to project a vector onto the feasible set.
This amounts to solving a problem of the form
\[
\min \|\m{x} - \m{a}\|^2 \mbox{ subject to } \m{g} (\m{x}) \le \m{0} .
\]
We solved this problem using BAA.
An iteration of BAA reduces to the solution of a problem
with the following structure:
\begin{equation}\label{Primal}
\min \|\m{x} - \m{a} \|^2
\mbox{ subject to } \|\m{x} - \m{c}_i \|^2 \le r_i^2, \quad
i = 1, 2, \ldots , m .
\end{equation}
As in \cite{Lin03}, we solve these problems by forming the dual problem
\[
\max_{\g{\lambda} \ge \m{0}} \;\; \min_{\m{x}\in \mathbb{R}^n}
~ \|\m{x} - \m{a}\|^2 + \sum_{i=1}^m \lambda_i
\left( \|\m{x}-\m{c}_i\|^2 - r_i^2 \right) .
\]
After carrying out the inner minimization, this reduces to
\begin{equation}\label{Dual}
\max_{\g{\lambda} \ge \m{0}}
- \frac{\|\m{a} + \sum_{i=1}^m \lambda_i \m{c}_i\|^2}
{1 + \sum_{i=1}^m \lambda_i} + \sum_{i=1}^m \lambda_i
(\|\m{c}_i\|^2 - r_i^2) .
\end{equation}
If $\g{\lambda}$ solves the dual problem (\ref{Dual}),
then the associated solution of the primal problem (\ref{Primal}) is
\[
\m{x} = \frac{\m{a} + \sum_{i=1}^m \lambda_i \m{c}_i}
{1 + \sum_{i=1}^m \lambda_i} .
\]
Again, the dual problem (\ref{Dual}) is solved using the
active set algorithm (ASA) of \cite{hz05a}.

The test problems used in Tables \ref{tab2} and \ref{tab3}
were generated as follows:
\def\baselinestretch{1}
\begin{table}[!t]
\caption{Positive definite cases }
\begin{small}
\begin{center}
\begin{tabular}
{|c|c|c|c|c|c|c|c|c|c|}

\hline $n,m$& & SED & success & BAA & success & NGPA & success &SPG & success \\
\hline  & \textit{time} & 0.52 & & 0.07 & &4.70 & & 5.94 & \\
  100,4 & \textit{iter} & 10.06 & 28 & 19.00 & 30  & 172.43 & 17 & 200.06 & 17 \\

\hline  & \textit{time} & 2.75 & & 0.32 & &10.68 & & 11.76 & \\
  200,4 & \textit{iter} & 9.56 & 26 & 12.70 & 30  & 203.23 & 21 & 348.36 & 20 \\

\hline  & \textit{time} & 8.14 & & 0.72 & &128.19 & & 122.23 & \\
  300,4 & \textit{iter} & 9.60 & 27 & 21.46 & 30  & 269.86 & 20 & 431.83 & 20 \\

\hline  & \textit{time} & 20.13 & & 1.64 & &404.28 & & 438.84 & \\
  400,4 & \textit{iter} & 10.26 & 27 & 49.26 & 29 & 352.66 & 18 & 545.23 & 18 \\

\hline  & \textit{time} & 44.07 & & 2.54 & &579.21 & & 647.56 & \\
  500,4 & \textit{iter} & 12.33 & 28 & 29.90 & 30  & 369.20 & 15 & 574.13 & 13 \\

\hline  & \textit{time} & 57.28 & & 4.27 & & 648.79 & & 611.60 & \\
  600,4 & \textit{iter} & 9.30 & 26 & 36.80 & 29  & 309.43 & 19 & 306.50 & 19 \\

\hline  & \textit{time} & 3.51 & & 0.08 & & 86.89 & & 81.67 & \\
 100,40 & \textit{iter} & 10.26 & 28 & 19.00 & 30 & 150.76 & 21 & 165.66 & 21 \\

\hline  & \textit{time} & 26.56 & & 0.32 & & 268.63 & & 250.22 & \\
 200,40 & \textit{iter} & 12.70 & 30 & 12.70 & 30  & 199.70 & 17 & 218.50 & 16 \\

\hline  & \textit{time} & 54.56 & & 0.81 & & 732.50 & & 727.92 & \\
200,100 & \textit{iter} & 10.66 & 30 & 9.50 & 30  & 295.80 & 20 & 327.26 & 20 \\

\hline  & \textit{time} & 23.84 & & 0.72 & & 579.62 & & 530.02 & \\
100,200 & \textit{iter} & 14.96 & 30 & 20.06 & 30  & 249.43 & 18 & 261.46 & 19 \\

\hline  & \textit{time} & 0.093 & & 0.002 & & 3.02 & & 2.75 & \\
  4,100 & \textit{iter} & 9.06 & 29 & 6.96 & 30  & 19.46 & 26 & 19.40 & 25 \\

\hline  & \textit{time} & 0.114 & & 0.004 & & 6.23 & & 5.70 & \\
  4,200 & \textit{iter} & 9.56 & 27 & 8.73 & 30  & 16.26 & 26 & 16.46 & 26 \\

\hline  & \textit{time} & 0.148 & & 0.012 & & 13.45 & & 11.58 & \\
  4,300 & \textit{iter} & 11.06 & 25 & 12.56 & 30  & 15.26 & 24 & 15.33 & 23 \\

\hline  & \textit{time} & 0.195 & & 0.014 & & 16.27 & & 12.87 & \\
  4,400 & \textit{iter} & 13.33 & 28 & 12.26 & 30  & 16.26 & 28 & 15.70 & 28 \\

\hline  & \textit{time} & 0.221 & & 0.017 & & 21.08 & & 18.16 & \\
  4,500 & \textit{iter} & 13.83 & 26 & 11.50 & 30  & 13.83 & 26 & 13.83 & 26 \\

\hline  & \textit{time} & 0.235 & & 0.018 & & 31.65 & & 34.83 & \\
  4,600 & \textit{iter} & 12.13 & 26 & 11.00 & 30  & 15.40 & 24 & 16.33 & 24 \\

\hline
\end{tabular}
\label{tab2}
\end{center}
\end{small}
\end{table}
\def\baselinestretch{1}
Let Rand$(n,l,u)$ denote a vector in $\mathbb{R}^n$
whose entries are chosen randomly in the interval $(l,u)$.
Random positive definite matrices $\m{A}$
are generated using the procedure given in
\cite{LinHan03}, which we now summarize.
Let $\m{w}_i \in$ Rand$(n, -1, 1)$ for $i = 1$, 2, 3, and define
\[
\m{Q}_i = \m{I} - 2\m{v}_i\m{v}_i\tr,
\quad \m{v}_i = \m{w}_i/\|\m{w}_i\| .
\]
Let $\m{D}$ be a diagonal matrix with diagonal in Rand$(n, 0, 100)$.
Finally, $\m{A} = \m{UDU}\tr$ with $\m{U} = \m{Q}_1\m{Q}_2\m{Q}_3$.
To obtain a randomly generated positive semidefinite matrix,
we use the same procedure, however, we randomly set one diagonal
element of $\m{D}$ to zero.

We make a special choice for $c_i$ to ensure that
the feasible set for (QP) is nonempty.
We first generate $\m{p} \in$ Rand$(n, -50, 50)$ and we set
\[
c_i = -(\m{p}\tr\m{A}_i\m{p} + \m{b}_i\tr\m{p} + s_i),
\]
where $s_i$ is randomly generated in the interval $[0, 10]$
and $\m{b}_i \in$ Rand$(n, -100, 100)$.
With this choice for $c_i$, the feasible set for (QP)
is nonempty since $\m{p}$ lies in the interior of the feasible set.
The stopping criterion in our experiments was
\begin{equation}\label{tolerance}
\|P(\m{x}_k - \m{g}_k) - \m{x}_k \| \le 10^{-4},
\end{equation}
where $P$ denotes projection into the feasible set for (QP) and
$\m{g}_k = 2\m{A}_0\m{x}_k + \m{b}_0$ is the gradient of the
objective function at $\m{x}_k$.
When the cost is convex,
the left side of (\ref{tolerance}) vanishes if and only if
$\m{x}_k$ is a solution of (QP).

Tables \ref{tab2} and \ref{tab3} report the average CPU time in
seconds ($time$), the average number of iterations ($iter$),
and the number of successes in 30 randomly generated test problems.
The algorithm was considered successful if the error tolerance
(\ref{tolerance}) was satisfied.
\def\baselinestretch{1}
\begin{table}[!t]
\caption{Positive semidefinite cases }

\begin{small}
\begin{center}
\begin{tabular}
{|c|c|c|c|c|c|c|c|}
\hline $n,m$& & BAA &success& NGPA &success& SPG &success\\
\hline  & \textit{time} & 0.11 & & 15.04 && 17.38 &\\
  100,4 & \textit{iter} & 42.16 & 30 & 328.53 & 22 & 408.20 & 19\\

\hline  & \textit{time} & 0.55 & & 44.32 && 45.82 &\\
  200,4 & \textit{iter} & 99.33 & 30 & 313.10 & 20 & 356.43 & 22\\

\hline  & \textit{time} & 1.02 & & 290.19 && 304.89 &\\
  300,4 & \textit{iter} & 74.23 & 30 & 374.13 & 22 & 417.60 & 21\\

\hline  & \textit{time} & 2.28 & & 501.14 && 572.63 &\\
  400,4 & \textit{iter} & 111.83 & 30 & 404.66 & 19 & 492.83 & 19\\

\hline  & \textit{time} & 5.37 & & 620.13 && 657.61 &\\
  500,4 & \textit{iter} & 200.03 & 27 & 382.66 & 16 & 478.60 & 17\\

\hline  & \textit{time} & 0.61 & & 356.40 && 321.92 &\\
 100,40 & \textit{iter} & 82.30 & 30 & 276.56 & 22 & 237.23 & 21\\

\hline  & \textit{time} & 2.74 & & 398.54 && 415.28 &\\
 200,40 & \textit{iter} & 127.23 & 30 & 369.43 & 17 & 416.060 & 17\\

\hline  & \textit{time} & 3.19 & & 1030.02 && 949.29 &\\
100,200 & \textit{iter} & 108.63 & 28 & 311.40 & 19 & 352.23 & 18\\

\hline  & \textit{time} & 0.054 & & 16.74 && 14.25 &\\
  4,100 & \textit{iter} & 38.66 & 30 & 31.13 & 16 & 31.23 & 16\\

\hline  & \textit{time} & 0.075 & & 44.74 && 33.24 &\\
  4,200 & \textit{iter} & 43.46 & 27 & 26.20 & 13 & 23.36 & 18\\

\hline  & \textit{time} & 0.076 & & 111.03 && 100.38 &\\
  4,300 & \textit{iter} & 31.50 & 29 & 29.33 & 14 & 28.60 & 11\\

\hline  & \textit{time} & 0.049 & & 205.17 && 237.62 &\\
  4,400 & \textit{iter} & 36.73 & 29 & 27.20 & 18 & 31.23 & 18\\

\hline  & \textit{time} & 0.065 & & 229.77 && 247.82 &\\
  4,500 & \textit{iter} & 41.86 & 28 & 24.46 & 16 & 26.30 & 17\\

\hline

\end{tabular}
\label{tab3}
\end{center}
\end{small}
\end{table}
\def\baselinestretch{1}

Based on our numerical experiments, it appears that
BAA can achieve an error tolerance on the order
of the square root of the machine epsilon \cite{Hager88, trefethen},
similar to the computing precision which is achieved by interior
point methods for linear programming prior to simplex crossover.
The convergence tolerance (\ref{tolerance}) was chosen since it
seems to approach the maximum accuracy which could be achieved by BAA
in these test problems.
Numerically, BAA seems to terminate when the
solution to the subproblem (\ref{eqa2}) yields a direction
which departs from the feasible set, and hence, the stepsize in the
line search Step 2b is zero.
We were able to achieve a further improvement in the solution
by taking a partial step in this infeasible direction
since the increase in constraint violation was much less than the
improvement in objective function value.
Nonetheless, the improvement in accuracy achieved by permitting infeasibility
was at most one digit in our experiments.

In Tables \ref{tab2} and \ref{tab3} we see that BAA gave the best results
for this test set, both in terms of CPU time and in terms of
successes (the number of times that the convergence tolerance
(\ref{tolerance}) was achieved).
Recall that the gradient projection algorithms in our experiments
used BAA to compute the projected gradient.
The convergence failures for the gradient projection algorithms in
Tables \ref{tab2} and \ref{tab3} were due to the fact that
BAA was unable to compute the projected gradient with enough
accuracy to yield descent in the gradient projection algorithm.

\subsection{Problems with nonconvex cost}
We tested our ellipsoidal branch and bound algorithm using some
randomly generated test problems with $\m{A}_0$ indefinite. To
compute $\g{\mu}$ in (\ref{mu}), we used the power method (see
\cite{trefethen}) to find the eigenvector associated with the
largest eigenvalue. We chose $\lambda$ in (\ref{fL}) to be 0.1 minus
the smallest eigenvalue of $\m{A}_0$. NGPA  was used to locally
solve (QP) and update the upper bound.

We took $m = 2$ and randomly generated test problem
using the procedure in \cite{Le00}.
That is, the ellipsoidal constraint functions in (QP) have
the form
\[
g_i (\m{x}) = (\m{x}-\m{c}_i)\tr\m{B}^{-1}_i(\m{x}-\m{c}_i) - 1,
\]
where $\m{B}_i = \m{UD}_i\m{U}\tr$ and $\m{U}$ is as given earlier.
$\m{D}_i$ is a diagonal matrix with its diagonal in Rand$(n, 0,
60)$, $\m{c}_1 \in$ Rand$(n,0,100)$, and $\m{c}_2 = \m{c}_1 +
.8\m{v}$ where $\m{v}$ is the semi-major axis of the ellipsoid $g_1
(\m{x}) \le 0$. For this choice of $\m{c}_2$, the ellipsoids $g_1
(\m{x}) \le 0$ and $g_2 (\m{x}) \le 0$ have nonempty intersection at
$\m{x} = \m{c}_2$. In the objective function, $\m{A}_0 =
\m{U}\m{D}\m{U}\tr$ where $\m{D}$ is a diagonal matrix with diagonal
in Rand$(n,-30,30)$ and $\m{b}_0 \in$ Rand$(n,-1,1)$. The case $m =
2$ is especially important since quadratic problems with two
ellipsoidal constraints belong to the class of Celis-Dennis-Tapia
subproblems \cite{cdt85} which arise from the application of the
trust region method for equality constrained nonlinear programming
\cite{Hei94, Yuan97, Yuan99a, Yuan99b, Yuan2005, Ye2003,
Martinez1995}.


If $\mbox{UB}_k$ and $\mbox{LB}_k$ are the respective upper and
lower bounds for the optimal objective function value at
iteration $k$, then our stopping criterion was
\[
\mbox{UB}_k - \mbox{LB}_k \le
\max \{\epsilon_a , \epsilon_r |\mbox{LB}_k|\},
\]
with $\epsilon_a = 10^{-5}$ and $\epsilon_r = 10^{-2}$.

We considered problems of 8 different dimensions ranging from 30 up
to 300 as shown in Table \ref{tab4}. For each dimension, we solved 4
randomly generated problems. Table \ref{tab4} shows the numerical
results for our test instances, where ``$neigs$'' is the number of
negative eigenvalues of the objective function, ``$lb^1$'' and
``$ub^1$'' are the lower bound and upper bounds at the first step,
``\textit{val}'' is the computed optimal value and ``\textit{it}''
is the number of iterations. We also report the performance of the
algorithm for $m = 6$ in Table \ref{tab5}.
\def\baselinestretch{1.0}
\begin{table}[!h]
\caption{The performance of branch and bound algorithm for $m = 2$}

\begin{small}
\begin{center}
\begin{tabular}
{|r|r|r|r|r|r|r|}
\hline $n$& \textit{$neigs$} & \textit{$lb^1$}& \textit{$ub^1$}&
\textit{val}&
\textit{it} & \textit{time} \\

\hline \hline
30 & 12 & 34827.3 & 35256.3 & 35254.8 & 5 & 1.75 \\

\hline
& 17 & -41212.1& -40746.2 & -40748.8 & 21 & 3.58 \\

\hline & 14 & 38601.4 & 38977.6 & 38977.6 & 0 & 0.72 \\

\hline
& 17 & -31534.2& -31108.4 & -31119.8 & 92 & 8.82 \\

\hline  50 & 22 &-357168.8 & -356828.8 & -356828.8 & 0 & 0.52 \\

\hline & 21 &-33792.9 & -33447.9 & -33447.9 & 1 & 0.84 \\

\hline & 21 & -29694.6 & -29254.1 & -29255.2 &  247 &  23.08 \\

\hline & 26 & 35034.0 & 35416.6 & 35414.8 & 5 & 2.12 \\

\hline 60 & 29 & 17783.6 & 18227.4 & 18227.4& 78 & 22.41 \\

\hline & 26 & -27498.2 & -27110.5 & -27110.5 & 69 & 18.22 \\

\hline &30 & -69845.7&-69463.1&-69463.1& 0 & 0.56 \\

\hline & 28& 20408.7 & 20963.1 &  20927.1 & 273 & 42.65 \\

\hline 100 & 50 & -11495.2 & -11196.9 & -11218.9    &   56 & 30.72 \\

\hline & 51 & 17539.6 &   17909.3 &   17909.3 &   4   &   1.84 \\

\hline &  52 &-46065.5  &  -45653.2  &  -45653.2  &  0 & 0.88 \\

\hline & 40& 970829.8    &   971326.6    &   971326.6    &   0   & 0.92 \\

\hline 150 & 75 & -302382.2   &   -302071.0   &   -302071.0   &   0   & 0.95 \\

\hline & 83& 29089.4 &   29500.8 &   29500.8 &   64  &   31.45 \\

\hline & 72& 16580.5 &   16904.9 &   16904.9 &   1   &    1.98  \\

\hline & 73& -32461.9    &   -32036.1    &   -32036.1    &   1   & 1.37 \\

\hline 200  & 100& 10798.5 &   11226.1 &   11226.1 &   81 &   56.58 \\

\hline & 95& -27242.9    &   -26792.1    &   -26792.1    &   2  & 2.27 \\

\hline & 100& 35293.0 &   35862.1 &   35862.1 &   1   &  1.63  \\

\hline & 96& -31712.8    &   -31138.3    &   -31138.3    &   77   & 47.06 \\

\hline 250 & 135 & 37015.8 &   37477.6 &   37477.6 &   1   &  2.9 \\

\hline & 131 & -27278.9    &   -26563.0   &   -26780.0   &   86 & 88.40 \\

\hline & 128 & -9979.6 &   -9683.9 &   -9683.9 &   59&  131.54 \\

\hline & 121 & -371385.9   &   -370991.9   &   -370991.9   &   0   & 2.03 \\

\hline 300 & 145 & -162041.5   &   -161645.7   &   -161645.7   &   0 & 5.33 \\

\hline & 152& -48085.4    &   -47529.3    &   -47529.3    &   1  & 7.56 \\

\hline & 138& 226345.6    &   226377.8    &   226377.8    &   0  & 4.79 \\

\hline & 148& -17649.5    &   -17013.7    &   -17323.2    &   109& 257.52 \\

\hline
\end{tabular}
\label{tab4}
\end{center}
\end{small}
\end{table}
\def\baselinestretch{1}

In comparing our ellipsoidal branch and bound algorithm based on
linear underestimation (EBL) to the ellipsoidal branch and bound
algorithm of Le Thi Hoai An \cite{Le00} based on dual
underestimation (EBD), an advantage of EBD is that the
underestimates are often quite tight in the dual-based approach. As
seen in Table \ref{tab4}, EBL required up to 273 bisections for this
test set while EBD in \cite{Le00} was able to solve randomly
generated test problems without any bisections. On the other hand, a
disadvantage of EBD is that the dual problems are nondifferentiable
when $\m{A}_0$ is indefinite. Consequently, the evaluation of the
lower bound using EBD entails solving an optimization problem which,
in general, is nondifferentiable. With EBL, however, computing a
lower bound involves solving a convex optimization problem. To
summarize, EBD provides tight lower bounds using a nondifferentiable
optimization problem for the lower bound, while EBL provides less
tight lower bounds using a convex optimization problem for the lower
bound.

\def\baselinestretch{1.0}
\begin{table}[!h]
\caption{The performance of branch and bound algorithm for $m = 6$}

\begin{small}
\begin{center}
\begin{tabular}
{|r|r|r|r|r|r|r|}

\hline $n$& \textit{$neigs$} & \textit{$lb^1$}& \textit{$ub^1$}&
\textit{val}&
\textit{it} & \textit{time} \\

\hline
\hline
30 & 17 & 22717.1 & 22993.1 & 22993.1 & 1 & 2.41 \\

\hline
& 21 & -22847.0& -22520.4 & -22524.5 & 14 & 10.58 \\

\hline & 20 & -17858.2 & -17573.2 & -17573.2 & 1 & 1.84 \\

\hline
60 & 33 & -21818.1 & -21489.1 & -21489.2 & 33 & 21.64 \\

\hline & 27 &47683.8 & 47826.0 & 47826.0 & 0 & 2.82 \\

\hline & 31 &-4926.5 & -4652.0 & -4728.7 & 4 & 7.62 \\

\hline
100 & 56 & -35438.9 & -35411.0 & -35411.0 & 0 & 0.78 \\

\hline & 52 & -1740.1 & -1187.5 & -1198.2 & 354 & 283.25 \\

\hline & 49 & -6756.5 &-6148.9&-6148.9& 3 & 8.06 \\

\hline
\end{tabular}
\label{tab5}
\end{center}
\end{small}
\end{table}
\def\baselinestretch{1}

\section{Conclusions}
\label{conclusions}
A globally convergent branch and bound algorithm was developed in
which the objective function was written as the difference of convex
functions. The algorithm was based on an affine underestimate given
in Theorem \ref{UnderTheorem} for the concave part of the objective
function restricted to an ellipsoid. An algorithm of Lin and Han
\cite{Lin03,LinHan03} for projecting a point onto a convex set was
generalized so as to replace their norm objective by an arbitrary
convex function. This generalization could be employed in the branch
and bound algorithm for a general objective function when the
constraints are convex. Numerical experiments were given for a
randomly generated quadratic objective function and randomly
generated convex, quadratic constraints.

\end{document}